\documentclass[final,3p,lefttitle,fleqn,12pt]{elsarticle}

\usepackage{amsthm}
\usepackage{amsmath}
\usepackage{amssymb}
\usepackage{mathrsfs}
\usepackage{mathdots}
\usepackage{latexsym}
\usepackage{amsfonts}
\usepackage{microtype}
\usepackage[colorlinks]{hyperref}
\usepackage{ulem}
\usepackage{multirow}
\usepackage{hhline}

\usepackage{xcolor}

 \usepackage{algorithm}
 \usepackage{algorithmic}
\usepackage[algo2e,ruled,linesnumbered]{algorithm2e}
\newcommand{\nonl}{%
  \let\algSavedNl\nl
  \renewcommand{\nl}{\let\nl\algSavedNl}%
}

\numberwithin{equation}{section}
\theoremstyle{plain}
\newtheorem{theorem}{Theorem}[section]
\newtheorem*{theorem*}{Theorem}

\newtheorem{lemma}[theorem]{Lemma}

\theoremstyle{remark}

\newtheorem{remark}[theorem]{\bf Remark}

\def\C{\mathbb C}

\def\det{\mathrm{det}\,}
\def\gcd{\mathrm{gcd}}

\def\0{\mathbf{0}}
\def\T{\mathrm{T}}
\def\A{\mathscr{A}}
\def\P{\mathscr{P}}

\def\p{\mathbf{p}}

\newcommand{\Ker}{\operatorname{Ker}}
\newcommand{\rank}{\operatorname{rank}}

\allowdisplaybreaks
\begin{document}

\begin{frontmatter}
\title{A new formula for Moore-Penrose inverses of rank-deficient Toeplitz matrices}

\author[label1]{Bohui Ban}
\ead{bhban@bnu.edu.cn}

\affiliation[label1]{organization={School of Mathematical Sciences, Laboratory of Mathematics and Complex System, MOE},
           addressline={Beijing Normal University},
             city={Beijing},
            postcode={100875},
            country={P.R. China}}

\affiliation[label2]{organization={Department of Mathematics, Faculty of Arts and Sciences},
         addressline={Beijing Normal University at Zhuhai},
         city={Zhuhai},
            postcode={519087},
            country={P.R. China}}

\author[label1,label2]{Yongjian Hu\corref{cor1}}
\ead{yongjian@bnu.edu.cn}

\cortext[cor1]{Corresponding author.}

\pdfstringdefDisableCommands{%
  \def\corref#1{<#1>}%
}

\date{\today}

\begin{abstract}
In this paper, we introduce a new formula for the Moore-Penrose inverses of the rank-deficient Toeplitz matrices, based on the quasi-direct decompositions of Toeplitz matrices. Numerical examples are presented to verify the validity and the efficiency of our method.
\end{abstract}



\begin{keyword} Toeplitz matrix \sep  Moore-Penrose inverse \sep Quasi-direct decomposition \sep Characteristic degree\sep Characteristic polynomial\\[1mm]
\MSC[2020]  15B05 \sep 15A09  \sep 

\end{keyword}

\end{frontmatter}


\section{Introduction}\label{Section1}

Let $T=(a_{i-j})_{i,j=0}^{m-1,n-1}\in\C^{m\times n}$ be a Toeplitz matrix. The unique solution $X\in\C^{n\times m}$ of the following four equations:
\begin{align}
     TXT=T,\quad XTX=X,\quad (TX)^*=TX,\quad (XT)^*=XT,\label{eq1.1}
\end{align}
is called the Moore-Penrose inverse of $T$, denoted by $T^{\dag}$. 
Any matrix satisfying the first two equations in \eqref{eq1.1} is called a $(1,2)$-inverse of $T$. And, the set of all $(1,2)$-inverses of $T$ is denoted by $T\{1,2\}$.

As the core tool for solving Toeplitz least squares problems, the Moore-Penrose inverse of $T$ has been applied in many fields, such as image restoration\cite{NP2014,FB1999}. 

The inverse of a nonsingular Toeplitz matrix is usually not a Toeplitz matrix. However, it can be represented as the sum of the products of lower triangular Toeplitz matrices and upper triangular Toeplitz matrices (see, e.g., \cite{HR1984,GS1972}). 
For a singular square Toeplitz matrix $T$, Xu \cite{X1992} gave a formula for the Moore-Penrose inverse based on the formula for the inverse of a nonsingular Toeplitz matrix. 
This formula, referred to as the Gohberg-Semencul-Trench type formula, utilizes the last column and last row of $T^{\dag}$, as well as the Moore-Penrose inverse of the $(n-1)$th leading principal submatrix of $T$, to represent the Moore-Penrose inverse as the sum of the products of lower triangular Toeplitz matrices and upper triangular Toeplitz matrices. 

Another formula for the Moore-Penrose inverse of a Toeplitz matrix is given via the Bezoutian, by generalizing the Bezoutian representations of the inverses of the Toeplitz and the Hankel matrices \cite{HR2010}. 
In 1994, Heinig and Hellinger \cite{HH1994} expressed the Moore-Penrose inverse of a square Toeplitz matrix as the sum of two Toeplitz-Bezoutians, which provided the possibility for developing fast algorithms to compute the Moore-Penrose inverse of Toeplitz matrices. 

Then, Heinig \cite{H2004}, using the augmented matrix approach, presented a formula for the Moore-Penrose inverse of a full column rank Toeplitz matrix $T$, based on the canonical fundamental system. With the help of this formula, Heinig gave new direct algorithms for the solution of $m \times n $ Toeplitz least squares problems with computational complexity $\mathcal{O}(n^2)+\mathcal{O}(m\mathrm{log}m)$. 

As for the calculation of the Moore-Penrose inverses of the Toeplitz matrices, Bini et al. \cite{BCB2003} utilized the Newton iteration method to compute the Moore-Penrose inverse of an $m \times n$ full-rank Toeplitz matrix in 2003. Subsequently, Wei et al. \cite{WCN2004} improved the algorithm proposed in \cite{BCB2003} for computing the Moore-Penrose inverse of an $m \times n$ rank-deficient Toeplitz matrix.
In 2016, Xie and Wei \cite{XW2016} presented a stability analysis of the Gohberg-Semencul-Trench type formula \cite{X1992} for the Moore-Penrose inverse of singular Toeplitz matrices and developed a fast algorithm for computing the Moore-Penrose inverse of Toeplitz matrices based on this formula. Compared to the Newton iteration method \cite{WCN2004}, this algorithm significantly improves computational efficiency when calculating the Moore-Penrose inverse of Toeplitz matrices, while maintaining the same level of output accuracy as the Newton iteration method.

In this paper, we first investigate the quasi-direct decomposition of rank-deficient Toeplitz matrices, using the theory of rational generating functions for Toeplitz matrices and prove the uniqueness of this decomposition, see Theorem\ref{thm.QDdem}. Using this decomposition, we then derive a new formula for the Moore-Penrose inverse of a rank-deficient Toeplitz matrix, see Theorem~\ref{thm.MPinv}.

The proposed formula applies to all rank-deficient Toeplitz matrices over the complex field. Compared with existing methods for computing the Moore-Penrose inverse of a Toeplitz matrix, this formula requires only the coefficient vector of the first characteristic polynomial and selected entries of the matrix itself. 

Furthermore, the formula provides a basis for designing an algorithm to compute the Moore-Penrose inverse of a rank-deficient Toeplitz matrix. 
Using the rank-deficient Toeplitz matrices from Section 5 of \cite{XW2016} as test cases, we compare the algorithm based on Theorem~\ref{thm.MPinv} with the F{\footnotesize M}MP algorithm presented in the reference \cite{XW2016}. The results show that, although the proposed algorithm is slightly slower for larger matrices, its computational errors are approximately four orders of magnitude smaller.

The remainder of this paper is organized as follows. In Section \ref{Section2}, we introduce the quasi-direct decomposition of the  rank-deficient Toeplitz matrix, and the existence and uniqueness of this decomposition are proved. 
In Section \ref{Section3}, we  derive a new formula for the Moore-Penrose inverse of the rank-deficient Toeplitz matrices. 
Numerical examples are presented in Section \ref{Section4}. Finally, we end with a summary and conclusions in Section \ref{Section5}.

\section{Quasi-direct decompositions of rank-deficient Toeplitz matrices}\label{Section2}

In this section, we introduce the quasi-direct decompositions of rank-deficient Toeplitz matrices and prove the existence and uniqueness of this decomposition, based on the theory of the rational generating functions of Toeplitz matrices. 

We first recall the definitions of the characteristic degrees and the first characteristic polynomial of Toeplitz matrix, mainly given in \cite{HR1984,BZH2026}.

Let $T=(a_{i-j})_{i,j=0}^{m-1,n-1}\in\C^{m\times n}$ be a nonzero Toeplitz matrix. We define a $k\times l$ ($k+l=m+n$) Toeplitz matrix as follows: 
\begin{equation*}
    T_{k,l}=(a_{i-j+m-k})_{i,j=0}^{k-1,l-1}=\left(
     \begin{array}{cccc}
       a_{m-k} & a_{m-k-1} & \cdots & a_{1-n} \\[3pt]
       a_{m-k+1} & a_{m-k} &\cdots &a_{2-n}   \\[3pt]
       \vdots &  \vdots &  &  \vdots \\[3pt]
       a_{m-1} &  a_{m-2} & \cdots & a_{k-n} \\
     \end{array}
   \right).
\end{equation*}
We define the first and second characteristic degrees $n_1$ and $n_2$ of the Toeplitz matrix $T$ as follows:  
\begin{align*}
  n_1=\left\{
\begin{aligned}
&\mathrm{rank\,}T_{k,k}, &\quad& {\mathrm if}\ m+n=2k,\\
& \mathrm{rank\,}T_{k,k+1}, &\quad& {\mathrm if}\ m+n=2k+1, \\
\end{aligned}
\right. \qquad
 n_2=m+n-n_1.
\end{align*}
It is obvious that $0<n_1\leq n_2$.

It has been shown in \cite{HR1984,BZH2026} that 
\begin{equation}\label{eq2.1}
  \rank T_{k,l} =\mathrm{min\,}\{k,l,n_1\},\quad k+l=m+n.
\end{equation}

Then, for any rank-deficient Toeplitz matrix $T$, we have $\rank T=\min\{m,n,n_1\}<\min\{m,n\}$, which implies that the first characteristic degree $n_1$ equals to the rank of $T$, that is to say $n_1=\rank T$. 

Then, we introduce the first characteristic polynomial of the Toeplitz matrix $T$. Let
\begin{equation}\label{eq2.2}
\A_{l}=(1,z,\cdots,z^{l})\Ker T_{k-1,l+1}.
\end{equation}
It can easily be seen that $\dim \A_{n_1}=\dim \Ker T_{n_2-1,n_1+1}=n_1+1-\rank T_{n_2-1,n_1+1}=n_1+1-\min\{n_2-1,n_1+1,n_1\}=1$. In this case, there exists a nonzero polynomial $p(z)$ of degree at most $n_1$ forming a basis of $\A_{n_1}$, which is called the first characteristic polynomial of $T$.

Suppose that
\begin{align}
    p(z)=p_{s}z^{s}+p_{s-1}z^{s-1}+\cdots+p_{t}z^t\quad(p_s,p_t\neq0,\quad 0\le t\le s\le n_1). \label{eq2.3}
\end{align} 
We denote by 
\begin{align}
    \p=(\underbrace{0,\cdots,0}_{t}, p_t, \cdots, p_s, \underbrace{0,\cdots,0}_{n_1-s})^{\T}\in\C^{n_1+1},
\end{align}
and 
\begin{align}
    \p_0=( p_t, \cdots, p_s )^{\T}\in\C^{s-t+1}.
\end{align}
Then, by $p(z)\in\A_{n_1}$, we have that $T_{n_2-1,n_1+1}\,\p=0$. 
We denote by $T_0$ the submatrix lying in columns from $t+1$ to $s+1$ of the matrix $T_{n_2-1,n_1+1}$, that is 
\begin{align*}
    T_0=\begin{pmatrix}
	a_{n_{1}+1-n-t} &  \cdots & a_{n_{1}+1-n-s}\\[2mm]
	\vdots & & \vdots\\
	a_{m-1-t} & \cdots & a_{m-1-s}\\
	\end{pmatrix}\in\C^{(n_2-1)\times(s-t+1)}.
\end{align*}
Then 
\begin{equation}\label{eq2.5}
T_0 \,\p_0=0. 
\end{equation}
We claim that $\dim\Ker T_0=1$. Indeed, when $n_1=n_2$, the second characteristic polynomial $q(z)$ of $T$ belongs to $\P^0_{n_1}$. If $\dim\Ker T_0>1$, then
\[
    \dim\Ker T_{n_1-1,n_1+1}\geq 1+\dim\Ker T_0>2,
\]
which contradicts $\rank T_{n_1-1,n_1+1}=n_1-1$.

When $n_1<n_2$, if $\dim\Ker T_0>1$, then
\[
    \dim\Ker T_{n_2-1,n_1+1}\geq\dim\Ker T_0>1,
\]
and hence $\rank T_{n_2-1,n_1+1}<n_1$, contradicting $\rank T_{n_2-1,n_1+1}=n_1$. Therefore, $\rank T_0=s-t$.

We next present a decomposition of a Toeplitz matrix that is key to constructing its quasi-direct decomposition. 

\begin{theorem}\label{thm.QTP}
    Let $T=(a_{i-j})_{i,j=0}^{m-1,n-1}\in\C^{m\times n}$ be a rank-deficient Toeplitz matrix, $n_1$ and $p(z)$ be the first characteristic degree and the first characteristic polynomial of $T$, respectively. If 
    \begin{align*}
           p(z)=p_{s}z^{s}+p_{s-1}z^{s-1}+\cdots+p_{t}z^t\quad(p_s,p_t\neq0,\quad 0\le t\le s\le n_1), 
    \end{align*}
    then
     \begin{align*}
    Q^{\T}TP=&\left(\begin{array}{c|c|cc}
         &   & O_{(n_1-s)\times(n-n_1)} &B\\[1mm]
         \hline 
         & A & &  \\
         \hline 
           O_{(m-n_1)\times t}& &   &   \\
         C & & &  
    \end{array}\right),
\end{align*}
where $A=(a_{n_1-t-s+i-j})_{i,j=0}^{s-t-1}$ is an invertible Toeplitz matrix, $B\in\C^{(n_1-s)\times (n_1-s)}$ and $C\in\C^{t\times t}$ are invertible upper triangular Toeplitz matrix and invertible lower  triangular Toeplitz matrix, and 
 \begin{align*}
    P=\begin{pmatrix}
    p_{t} & & & & & & \\
    \vdots & \ddots & & O_{t\times (s-t)} & & &  \\
      p_{s} & & p_{t} & & & & \\
          & \ddots & \vdots & I_{s-t}  & p_{t} & & \\
        & & p_{s} & & \vdots & \ddots & \\
       & & & & p_{s} &  &p_{t} \\
      & & & O_{(n-s)\times (s-t)} & & \ddots & \vdots \\
         & & & & & & p_{s}
   \end{pmatrix}
        \in \C^{n\times n},
\end{align*}
\begin{align*}
    Q=\begin{pmatrix}
    p_{s} & & & & & & \\
    \vdots & \ddots & & O_{(n_1-s)\times (s-t)} & & &  \\
      p_{t} & & p_{s} & & & & \\
          & \ddots & \vdots & I_{s-t}  & p_{s} & & \\
        & & p_{t} & & \vdots & \ddots & \\
       & & & & p_{t} &  &p_{s} \\
      & & & O_{(m-n_1+t)\times (s-t)} & & \ddots & \vdots \\
         & & & & & & p_{t}
   \end{pmatrix}
        \in \C^{m\times m}.
\end{align*}
 are invertible matrices. 
\end{theorem}

\begin{proof}
By \eqref{eq2.5}, we can calculate that
\begin{align*}
    Q^{\T}TP=&\setlength{\arraycolsep}{2pt}
    Q^{\T} \left(\begin{array}{ccc|ccc|cccc}
         & & &  a_{-t} & \cdots & a_{1-s} & & b'_{1} & \cdots & b'_{n_{1}-s}\\
         &  & & \vdots & & \vdots & O_{(n_1-s)\times(n-n_1)} & & \ddots & \vdots\\
         & & & a_{n_1-t-s-1} & \cdots & a_{n_1-2s} & & & & b'_{1}\\[1.5mm]
         \hline 
         & & & a_{n_1-t-s} & \cdots & a_{n_1-2s+1} & & & & \\
         &  & & \vdots & & \vdots &  & & & \\
         & & & a_{n_1-2t-1} & \cdots & a_{n_1-t-s} & & & & \\[1.5mm]
         \hline
          \multicolumn{3}{c|}{O_{(m-n_1)\times t}} & a_{n_1-2t} & \cdots & a_{n_1-t-s+1} & & & &\\
         c'_{1} & & &   &   &   & & & &\\
         \vdots & \ddots & & \vdots & & \vdots &  & & &\\
         c'_{t} & \cdots & c'_{1} & a_{m-t-1} & \cdots & a_{m-s-1} & & & &
    \end{array}\right)\\[0.3cm]
    =&\left(\begin{array}{ccc|c|cccc}
         & & & & & b_{1} & \cdots & b_{n_{1}-s}\\
         & & & & O_{(n_1-s)\times(n-n_1)} & & \ddots & \vdots\\
         & & & & & & & b_{1}\\[1mm]
         \hline 
         & & & A & & & & \\
         \hline
          \multicolumn{3}{c|}{O_{(m-n_1)\times t}} &    &   & & &\\
         c_{1} & & &   &   &   & & \\
         \vdots & \ddots & &  & &  &   &\\
         c_{t} & \cdots & c_{1} &   &  & & &
    \end{array}\right),
\end{align*}
where
\begin{align*}
    &\left(\begin{array}{c}
         b_1  \\
          \vdots\\
          b_{n_1-s}
    \end{array}\right)=\setlength{\arraycolsep}{1.8pt}\left(\begin{array}{ccc}
         p_s & &  \\
          \vdots&\ddots&\\
         p_{2s-n_1+1}&\cdots& p_s
    \end{array}\right)\left(\begin{array}{c}
         b'_1  \\
          \vdots\\
          b'_{n_1-s}
    \end{array}\right),\\ 
    & \qquad \qquad
    \left(\begin{array}{c}
         b'_1  \\
          \vdots\\
          b'_{n_1-s}
    \end{array}\right)=\setlength{\arraycolsep}{1.8pt}\left(\begin{array}{ccc}
         a_{n_1-n-t} &\cdots &a_{n_1-n-s}  \\
          \vdots& &\vdots\\
         a_{1-n-t+s}&\cdots& a_{1-n}
    \end{array}\right)\left(\begin{array}{c}
         p_t  \\
          \vdots\\
          p_s
    \end{array}\right),\\[2mm]
    &\left(\begin{array}{c}
         c_1  \\
          \vdots\\
          c_t
    \end{array}\right)=\setlength{\arraycolsep}{1.8pt}\left(\begin{array}{ccc}
         p_t & &  \\
          \vdots&\ddots&\\
         p_{2t-1}&\cdots& p_t
    \end{array}\right)\left(\begin{array}{c}
         c'_1  \\
          \vdots\\
          c'_t
    \end{array}\right),\\
   & \qquad \qquad
    \left(\begin{array}{c}
         c'_1  \\
          \vdots\\
          c'_t
    \end{array}\right)=\setlength{\arraycolsep}{1.8pt}\left(\begin{array}{ccc}
         a_{m-t} &\cdots &a_{m-s}  \\
          \vdots& &\vdots\\
         a_{m-1}&\cdots& a_{m-s+t-1}
    \end{array}\right)\left(\begin{array}{c}
         p_t  \\
          \vdots\\
          p_s
    \end{array}\right),
\end{align*}
with $p_j=0$ for $j<t$ and $j>s$.

Let
\begin{align*}
    B=\left(\begin{array}{ccc}
      b_{1} & \cdots & b_{n_{1}-s}\\
         & \ddots & \vdots\\
       & & b_{1}\\[1mm]
    \end{array}\right), \qquad
    C=\left(\begin{array}{ccc}
         c_{1} & &\\
         \vdots & \ddots & \\
         c_{t} & \cdots & c_{1} 
    \end{array}\right). 
\end{align*}
We assume that $b_1=0$, then $b'_1=0$. So, 
\begin{align*}
\begin{pmatrix}
  a_{n_{1}-n-t} &  \cdots & a_{n_{1}-n-s}\\[1mm]
   &  T_0 &  
\end{pmatrix}\begin{pmatrix}
	p_{t}\\[2mm]
	\vdots\\
	p_{s}\\
	\end{pmatrix}=0, 
\end{align*}
Then
\begin{equation*}
  (\underbrace{0,\cdots,0}_{t},\,p_{t},\cdots, p_{s},\,\underbrace{0,\cdots,0}_{n_{1}-s-1})^{\T}\in\Ker T_{n_{2},n_{1}}. 
\end{equation*}
This is a contradiction of the fact that $\dim\Ker T_{n_{2},n_{1}}=0$. Hence $b_1\neq 0$. Then, $B$ is  invertible. Similarly, $C$ is invertible.

Finally, it remains to prove that $A=(a_{n_1-t-s+i-j})_{i,j=0}^{s-t-1}$ is an invertible Toeplitz matrix. Similarly, by using elementary transformations on $T_0$, we have that
\begin{align*}
   \left(\begin{array}{cccccccc}
    p_{s} &\cdots &p_t & & & & &\\
     & \ddots & &\ddots & & & & \\
     & & p_{s} & \cdots & p_{t} & & &  \\
         \multicolumn{3}{c}{O_{(s-t)\times(n-1-s) }}  & I_{s-t}  &   \multicolumn{4}{c}{O_{(s-t)\times(m-n_1+t)}} \\
        & && p_{s} & \cdots&p_t  & & \\
      & && & \ddots & & \ddots &  \\
        & & & & & p_{s} & \cdots &p_{t} 
   \end{array}\right) T_0&\left(\begin{array}{c|c}
          & p_t  \\
         I_{s-t} & \vdots  \\
         &p_{s-1}\\
         \hline
         & p_s
    \end{array}\right)\\=&\left(\begin{array}{c|c}
        O_{(n-1-s)\times (s-t)} & 0  \\
        A & \vdots\\
        O_{(m-n_1+t)\times (s-t)} &0
    \end{array}\right).
\end{align*}
Therefore, $\rank A=\rank T_0=s-t$, $A$ is invertible. 
\end{proof}

\begin{remark}
     If $s=t$, $n_1=s$, or $t=0$, the corresponding block $A$, $B$, or $C$, respectively, is omitted.
\end{remark}

\vspace{6mm}
Let $A$, $B$, $C\in\C^{m\times n}$. We call the matrix $A$ a quasi-direct sum of the matrices $B$ and $C$ if and only if
\begin{equation*}
        A=B+C,\quad \rank A=\rank B+\rank C.
\end{equation*}
This decomposition is called the quasi-direct decomposition of $A$. 

The quasi-direct decomposition of the Hankel matrix has been studied in \cite{F1984}. Fiedler \cite{F1984} has proved that any square Hankel matrix $H$ can be expressed as a quasi-direct sum of a proper  and a degenerate Hankel matrix, and that this decomposition is unique if and only if $H$ is singular. 

Since Hankel matrices and Toeplitz matrices can be transformed into each other via multiplication by a flip matrix, Fiedler \cite{F1984} has pointed out that Toeplitz matrices admit similar decompositions. 
Nevertheless, the quasi-direct decomposition of Hankel matrices may not be suitable for Toeplitz matrices. 
Although Hankel matrices have similar structure to Toeplitz matrices, the quasi-direct decompositions of Hankel matrices cannot be directly extended to Toeplitz matrices. One important reason is that there is a significant difference between the rational generating functions of Hankel matrices and Toeplitz matrices, see e.g. \cite{HR2010,HJ1986,HJ1986_1}. The rational generating functions of Hankel matrices require a given expansion series only at infinity, while the rational generating functions of Toeplitz matrices require the given expansion series at both  the origin and infinity. 
On the other hand, according to \cite{HR1984}, the first characteristic degree of the Toeplitz matrix is divided into three indicators, rather than being divided into two parts like the Hankel matrix. This seems to imply that the Toeplitz matrix should be the quasi-direct sum of three Toeplitz matrices with smaller rank. 

In order to study the quasi-direct decompositions of Toeplitz matrices, we begin by introducing the definition of proper Toeplitz matrices and degenerate Toeplitz matrices. 

Let $T=(a_{i-j})_{i,j=0}^{m-1,n-1}\in\C^{m\times n}$ be a Toeplitz matrix. If the maximum order of the nonsingular leading  principle submatrix of $T$ is equal to $\rank T$, then $T$ is proper; If $a_{-q+1}=\cdots=a_0=\cdots=a_{p-1}=0$ and $p+q>\max\{m,n\}$, that is to say $T$ has the following form:
\begin{align*}
 \left(\begin{array}{ccc|c|ccc}
          &   &  &  & a_{-q} & \cdots & a_{1-n}\\
	    & & & & & \ddots & \vdots \\
	    & & & & & & a_{-q}\\[1.5mm]
           \hline 
          & &  & O_{(p+q-n)\times(p+q-m)}&  & &  \\[1.5mm]
           \hline 
          a_{p}& &  &  &  & & \\
         \vdots & \ddots & & & & &\\
        a_{m-1} & \cdots & a_{p}&  &  & & 
 \end{array}\right),
\end{align*}
then $T$ is degenerate. 

It is obvious that all leading principal submatrices of the degenerate Toeplitz matrix are singular. The zero matrix is degenerate. In fact,  all leading principal submatrices of the zero matrix are also zero matrices. We can thus conclude that the maximum order of the nonsingular leading principle submatrix of the zero matrix is equal to $0$. Therefore, the zero matrix can be classified as a proper Toeplitz matrix.

We first present a property of the rank-deficient Toeplitz matrices.

\begin{lemma}\label{lem2.s-t}
     Let $T=(a_{i-j})_{i,j=0}^{m-1,n-1}\in\C^{m\times n}$ be a rank-deficient Toeplitz matrix, $n_1$ and $p(z)$ be the first characteristic degree and the first characteristic polynomial of $T$, respectively. If $\deg p(z)=s$, $\gcd (p(z),z^{n_1})=z^t$, then the maximum order of the non-singular leading principal submatrix of $T$ is $s-t$.
\end{lemma}

\begin{proof}
Let $T_k$ ($1\le k\le\min \{m,n\}$) be the leading principle submatrix of $T$ of order $k$. By using elementary transformations on $T_0$, we notice that
\begin{align*}
   \left(\begin{array}{ccccc}
         \multicolumn{2}{c}{O_{(s-t)\times(n+t-n_1-1) }}  & I_{s-t}  &   \multicolumn{2}{c}{O_{(s-t)\times(m-s)}} \\
       p_{s} & \cdots&p_t  & & \\
       & \ddots & & \ddots &  \\
         & & p_{s} & \cdots &p_{t} 
   \end{array}\right) T_0\left(\begin{array}{c|c}
          & p_t  \\
         I_{s-t} & \vdots  \\
         &p_{s-1}\\
         \hline
         & p_s
    \end{array}\right)\\
    =\left(\begin{array}{c|c}
     &0\\
        T_{s-t} & \vdots\\
        O_{(m-n_1+t)\times (s-t)} &0
    \end{array}\right).\label{eq2.6}
\end{align*}
Then,  $\rank T_{s-t}=\rank T_0=s-t$, and then $T_{s-t}$ is invertible. 

Then, we obtain that for $s-t+1 \leq j \leq m-t$, 
$$ 
T_j (p_t, \cdots, p_s, \underbrace{0,\cdots,0}_{j-(s-t+1)} )^{\T}=0,
$$
and for $m-t+1 \leq j \leq \min \{m, n\}$,  
$$ 
T_j (\underbrace{0,\cdots,0}_{j-(m-t)}, p_t, \cdots, p_s, \underbrace{0,\cdots,0}_{m-s-1})^{\T}=0.
$$ 
Hence, $\det T_j=0$ for $s-t+1 \leq j \leq\min\{m, n\}$.

\end{proof}

Next, we prove that any rank-deficient Toeplitz matrix can be expressed as a quasi-direct sum of a proper Toeplitz matrix and a degenerate Toeplitz matrix, and this decomposition is unique. 

\begin{theorem}\label{thm.QDdem}
    Let $T=(a_{i-j})_{i,j=0}^{m-1,n-1}\in\C^{m\times n}$ be a rank-deficient Toeplitz matrix. Then $T$ has the following decomposition:
    \begin{equation*}
        T=T_p+T_d,\quad \rank T=\rank T_p+\rank T_d,
    \end{equation*}
    where $T_p$ is a proper Toeplitz matrix, and $T_d$ is a degenerate Toeplitz matrix. Moreover, this decomposition is unique. 
\end{theorem}

\begin{proof}
Let $n_1$, $p(z)$ be the first characteristic degree and the first characteristic polynomial of $T$, and $p(z)$ be given by \eqref{eq2.3}. 
Let $P$, $Q$ be as defined in  Theorem \ref{thm.QTP}. 

Choose
\begin{equation*}
T_p=(Q^{\T})^{-1}\left(\begin{array}{c|c|c}
      &   &   O_{(n_1-s)\times (n-s)} \\[1mm]
         \hline
          & A &  \\[1mm]
         \hline
          O_{(m-n_1+t)\times t} & & 
    \end{array}\right)P^{-1}, 
\end{equation*}
where $A$ is same as in Theorem\ref{thm.QTP}. Then $T_p=(a^{'}_{i-j})_{i,j=0}^{m-1,n-1}$ is a Toeplitz matrix of rank $s-t$, where
 \begin{align*}
    &a_i^{'}=a_i, \quad i=n_1-s-n+1,\cdots,0,\cdots, m-t-1;\\
    &a_i^{'}=-\frac{1}{p_t}(a'_{i-1} p_{t+1}+\cdots+a'_{i+t-s}p_{s}), \quad i=m-t,\cdots, m-1;\\
    &a_i^{'}=-\frac{1}{p_s}(a'_{i+1} p_{s-1}+\cdots+a'_{i+s-t}p_{t}), \quad i=1-n,\cdots,n_1-s-n.
\end{align*}
In this case, the maximum order  of the nonsingular principle submatrix of $T_p$ is $s-t$. Therefore, $T_p$ is proper. 

Choose
\begin{equation*}
     T_d=(Q^{\T})^{-1}\left(\begin{array}{c|c|cc}
         &   & O_{(n_1-s)\times(n-n_1)} &B\\[1mm]
         \hline 
         & O_{s-t}& &  \\
         \hline 
           O_{(m-n_1)\times t}& &   &   \\
         C & & &  
    \end{array}\right)P^{-1},
\end{equation*}
in which $B$, $C$ are same as Theorem \ref{thm.QTP}.

Then
\begin{equation*}
    T_d=T-T_p=\left(\begin{array}{ccc|c|ccc}
             &  &  & & l_{1} & \cdots &l_{n_{1}-s}\\
            & &  &  & & \ddots & \vdots\\
            &  &  & & & &  l_{1}\\[1.5mm]
           \hline 
          & &  & O_{(m-n_{1}+s-t)\times(n-n_{1}+s-t)}&  &  &   \\[1.5mm]
           \hline 
         l'_1 & & &  &  &  & \\
         \vdots & \ddots & & & & & \\
         l'_{t} & \cdots & l'_{1} &  &  &  &  
 \end{array}\right)
\end{equation*}
is a degenerate Toeplitz matrix. We note that 
\begin{align*}
  &  l_1=a_{n_1-s-n}-a'_{n_1-s-n}=\frac{1}{p_s}b'_1\neq0,\\[1mm]
  &  l'_1=a_{m-t}-a'_{m-t}=\frac{1}{p_t}c'_1\neq0 ,
\end{align*}
Hence $\rank T_d=t+(n_1-s)$. Therefore, $T$ has decomposition $T=T_p+T_d$, where $T_p$ is proper and $T_d$ is degenerate, and $\rank T=n_1=\rank T_p+\rank T_d$.
 
Next, we begin to prove the uniqueness of the quasi-direct decomposition. Assume that the Toeplitz matrix $T$ has another decomposition $T=T^{'}_p+T^{'}_d$, in which $T^{'}_p$ is proper with rank $r_1$, and 
\begin{equation*}
    T^{'}_d=\left(\begin{array}{ccc|c|ccc}
             &  &  & & d_{1} & \cdots &d_{r_3}\\
            & &  &  & & \ddots & \vdots\\
            &  &  & & & &  d_{1}\\[1.5mm]
           \hline 
          & &  & O_{(m-r_2-r_3)\times(n-r_2-r_3)}&  &  &   \\[1.5mm]
           \hline 
         d'_1 & & &  &  &  & \\
         \vdots & \ddots & & & & & \\
         d'_{r_2} & \cdots & d'_{1} &  &  &  &  
 \end{array}\right)\ \ (d_1,d'_1\neq0) 
\end{equation*}
is degenerate with rank $r_2+r_3$. Then $r_1+r_2+r_3=\rank T=n_1<\min \{m , n\}$. So, $m-r_2> r_1+r_3 \geq r_1$, $n-r_3> r_1+r_2 \geq r_1$. 

Since $T^{'}_d$ is degenerate, $T$ and $T^{'}_p$ have the same leading principal submatrix, which is a proper Toeplitz matrix of order $\min\{m-r_2,n-r_3 \}$ with rank $r_1$. We note that 
\begin{align*}
    &\rank T_{\min\{m-r_2,n-r_3 \}+i}\leq r_1+i<\min\{m-r_2,n-r_3 \}+i,\\
   &\qquad\qquad\qquad i=1,\cdots,\max\{m-r_2,n-r_3 \}-\min\{m-r_2,n-r_3 \},
\end{align*}
and,      
\begin{align*}
    &\ \rank T_{\max\{m-r_2,n-r_3 \}+i}\\
    \leq&\ r_1+\max\{m-r_2,n-r_3 \}-\min\{m-r_2,n-r_3 \}+2i\\
    \leq &\ \max\{m-r_2,n-r_3 \}+i+r_1-\min\{m-r_2,n-r_3 \}+\\
    &\qquad\ \min\{m,n \}-\max\{m-r_2,n-r_3 \}\\
   = &\ \max\{m-r_2,n-r_3 \}+i+r_1-(m-r_2+n-r_3)+\min\{m,n \}\\
    = &\ \max\{m-r_2,n-r_3 \}+i+n_1-\max\{m,n \}\\
      < &\ \max\{m-r_2,n-r_3 \}+i,\\
    &\qquad\qquad\quad i=1,\cdots,\min\{m,n \}-\max\{m-r_2,n-r_3 \}. 
\end{align*}    
Therefore $\det T_j=0$, $r_1<j\leq \min\{m,n \}$. So, the maximum order of the nonsingular submatrices of $T$ is $r_1$. By Lemma \ref{lem2.s-t}, we have that $r_1=s-t$. We conclude that $T_p=T^{'}_p$, and hence that $T_d=T^{'}_d$. We prove the uniqueness of the quasi-direct decomposition of Toeplitz matrices. 
\end{proof}

\section{A formula for the Moore-Penrose inverse of rank-deficient Toeplitz matrices}\label{Section3}


In this section, we derive a new formula of the Moore-Penrose inverse of a nonzero rank-deficient Toeplitz matrix, utilizing the decomposition of Toeplitz matrices.


We first use Theorem~\ref{thm.QTP} to obtain a new formula for a $\{1,2\}$-inverse of a nonzero rank-deficient Toeplitz matrix.

\begin{theorem}\label{thm.T12}
    Let $T=(a_{i-j})_{i,j=0}^{m-1,n-1}\in\C^{m\times n}$ be a nonzero rank-deficient Toeplitz matrix. Then
\begin{align}
    P \left(\begin{array}{c|c|cc}
         &   &  O_{t\times (m-n_1)} &C^{-1}\\[1mm]
         \hline
         & A^{-1} & &  \\
         \hline
          O_{(n-n_1)\times(n_1-s)}& &   &   \\[1mm]
         B^{-1} & & &  
    \end{array}\right)Q^{\T} \in T\{1,2\},\label{T12}
\end{align}
where $A$, $B$, $C$, $P$, and $Q$ are as defined in Theorem~\ref{thm.QTP}.
\end{theorem}

\begin{proof}
    Let $X$ denote the matrix on the left-hand side of \eqref{T12}. Then
    \begin{align*}
        \setlength{\arraycolsep}{1pt}
        TXT=&(Q^{\T})^{-1}\setlength{\arraycolsep}{1.5pt}\left(\begin{array}{c|c|cc}
         &   & O_{(n_1-s)\times(n-n_1)} &B\\[1mm]
         \hline
         & A & &  \\
         \hline
           O_{(m-n_1)\times t}& &   &   \\
         C & & &  
    \end{array}\right)\\
   & \hspace{-1.1cm}\setlength{\arraycolsep}{0.5pt}\left(\begin{array}{c|c|cc}
         &   &  O_{t\times (m-n_1)} &C^{-1}\\[1mm]
         \hline
         & A^{-1} & &  \\
         \hline
          O_{(n-n_1)\times(n_1-s)}& &   &   \\
         B^{-1} & & &  
    \end{array}\right)\hspace{-2mm}\left(\begin{array}{c|c|cc}
         &   & O_{(n_1-s)\times(n-n_1)} &B\\
         \hline
         & A & &  \\
         \hline
           O_{(m-n_1)\times t}& &   &   \\
         C & & &  
    \end{array}\right)\hspace{-2mm}P^{-1}\\
    =&(Q^{\T})^{-1}\setlength{\arraycolsep}{1.5pt}\left(\begin{array}{c|c|cc}
         &   & O_{(n_1-s)\times(n-n_1)} &B\\
         \hline
         & A & &  \\
         \hline
           O_{(m-n_1)\times t}& &   &   \\
         C & & &  
    \end{array}\right) \left(\begin{array}{ccc}
      I_s   &   &  \\
         & O_{(n-n_1)}&   \\
         & & I_{n_1-s}  
    \end{array}\right)P^{-1}\\
     =&(Q^{\T})^{-1}\setlength{\arraycolsep}{1.5pt}\left(\begin{array}{c|c|cc}
         &   & O_{(n_1-s)\times(n-n_1)} &B\\
         \hline
         & A & &  \\
         \hline
           O_{(m-n_1)\times t}& &   &   \\
         C & & &  
    \end{array}\right) P^{-1}=T,
    \end{align*}
    and
    \begin{align*}
        \setlength{\arraycolsep}{1pt}
        XTX=&P \setlength{\arraycolsep}{1.5pt}\left(\begin{array}{c|c|cc}
         &   &  O_{t\times (m-n_1)} &C^{-1}\\
         \hline
         & A^{-1} & &  \\
         \hline
          O_{(n-n_1)\times(n_1-s)}& &   &   \\
         B^{-1} & & &  
    \end{array}\right)\\
   &\hspace{-1.1cm}\setlength{\arraycolsep}{0.5pt}\left(\begin{array}{c|c|cc}
         &   & O_{(n_1-s)\times(n-n_1)} &B\\
         \hline
         & A & &  \\
         \hline
           O_{(m-n_1)\times t}& &   &   \\
         C & & &  
    \end{array}\right)\hspace{-2mm}\left(\begin{array}{c|c|cc}
         &   &  O_{t\times (m-n_1)} &C^{-1}\\
         \hline
         & A^{-1} & &  \\
         \hline
          O_{(n-n_1)\times(n_1-s)}& &   &   \\
         B^{-1} & & &  
    \end{array}\right)\hspace{-1mm}Q^{\T}\\
    =&P\setlength{\arraycolsep}{1.5pt}\left(\begin{array}{ccc}
      I_s   &   &  \\
         & O_{(n-n_1)}&   \\
         & & I_{n_1-s}  
    \end{array}\right)\left(\begin{array}{c|c|cc}
         &   &  O_{t\times (m-n_1)} &C^{-1}\\
         \hline
         & A^{-1} & &  \\
         \hline
          O_{(n-n_1)\times(n_1-s)}& &   &   \\
         B^{-1} & & &  
    \end{array}\right) Q^{\T}\\
     =&P\setlength{\arraycolsep}{1.5pt}\left(\begin{array}{c|c|cc}
         &   &  O_{t\times (m-n_1)} &C^{-1}\\
         \hline
         & A^{-1} & &  \\
         \hline
          O_{(n-n_1)\times(n_1-s)}& &   &   \\
         B^{-1} & & &  
    \end{array}\right) Q^{\T}=X.
    \end{align*}
    Therefore, $X$ is a $\{1,2\}$-inverse of $T$, that is, $X\in T\{1,2\}$.
\end{proof}

\begin{remark}
    If $s=t$, $n_1=s$, or $t=0$, the corresponding block $A^{-1}$, $B^{-1}$, or $C^{-1}$, respectively, is omitted.
\end{remark}

The preceding result gives a formula for a $\{1,2\}$-inverse of $T$. We can then obtain its Moore-Penrose inverse by using orthogonal projections.

\begin{theorem}\label{thm.MPinv}
    Let $T=(a_{i-j})_{i,j=0}^{m-1,n-1}\in\C^{m\times n}$ be a nonzero rank-deficient Toeplitz matrix. Then
    \begin{align}
    T^{\dag}=(I_n-\widetilde{N})P \left(\begin{array}{c|c|cc}
         &   &  O_{t\times (m-n_1)} &C^{-1}\\[1mm]
         \hline
         & A^{-1} & &  \\
         \hline
          O_{(n-n_1)\times(n_1-s)}& &   &   \\[1mm]
         B^{-1} & & &  
    \end{array}\right)Q^{\T}(I_m-\widetilde{M}),\label{Tmp}
\end{align}
where $A$, $B$, $C$, $P$, and $Q$ are as defined in Theorem~\ref{thm.QTP}, and $\widetilde{M}\in\C^{m\times m}$ and $\widetilde{N}\in\C^{n\times n}$ are the orthogonal projection matrices onto $\Ker T^*$ and $\Ker T$, respectively.
\end{theorem}

\begin{proof}
Since $\widetilde{N}\in\C^{n\times n}$ is the orthogonal projection matrix onto $\Ker T$, for any vector $x\in\C^n$ we have $\widetilde{N}x\in\Ker T$ and $T\widetilde{N}x=T(\widetilde{N}x)=0$. Hence, $T\widetilde{N}=O$, which gives $T(I_n-\widetilde{N})=T$. Similarly, since $\widetilde{M}\in\C^{m\times m}$ is the orthogonal projection matrix onto $\Ker T^*$, we have $T^*(I_m-\widetilde{M})=T^*$, and therefore $(I_m-\widetilde{M})T=T$.

Let $X$ denote the matrix on the right-hand side of~\eqref{Tmp}, and let $G$ denote the matrix on the left-hand side of~\eqref{T12}. Then
\begin{align*}
    TXT=&T(I_n-\widetilde{N})G(I_m-\widetilde{M})T=TGT=T,\\
    XTX=&(I_n-\widetilde{N})G(I_m-\widetilde{M})T(I_n-\widetilde{N})G(I_m-\widetilde{M})\\
    =&(I_n-\widetilde{N})GT(I_n-\widetilde{N})G(I_m-\widetilde{M})\\
    =&(I_n-\widetilde{N})GTG(I_m-\widetilde{M})\\
    =&(I_n-\widetilde{N})G(I_m-\widetilde{M})=X.
\end{align*}

Since $TGT=T$, or equivalently $T(GT-I_n)=O$, we have $\Im (GT-I_n)\subseteq \Ker T$. Consequently, $\widetilde{N} (GT-I_n)=GT-I_n$, that is, $(I_n-\widetilde{N})GT=I_n-\widetilde{N}$. Thus,
\begin{align*}
    XT=(I_n-\widetilde{N})G(I_m-\widetilde{M})T=(I_n-\widetilde{N})GT=(I_n-\widetilde{N}), 
\end{align*}
which is Hermitian, that is, $(XT)^*=XT$.

Similarly, since $T^*(G^*T^*-I_m)=O$, we obtain $\widetilde{M} (G^*T^*-I_m)=G^*T^*-I_m$, that is, $(I_m-\widetilde{M})G^*T^*=I_m-\widetilde{M}$. Hence,
\begin{align*}
    (TX)^*=[T(I_n-\widetilde{N})G(I_m-\widetilde{M})]^*=[TG(I_m-\widetilde{M})]^*=I_m-\widetilde{M}, 
\end{align*}
which is Hermitian, and therefore $(TX)^*=TX$.
 
Thus, $X$ satisfies all four equations defining the Moore-Penrose inverse, and hence $X=T^{\dag}$.
\end{proof}

The preceding theorem provides a new formula for the Moore-Penrose inverse of a nonzero rank-deficient Toeplitz matrix. Formula~\eqref{Tmp} shows that $T^{\dag}$ can be obtained from a $\{1,2\}$-inverse by constructing two orthogonal projection matrices. 
To construct the orthogonal projection matrices $\widetilde{M}$ and $\widetilde{N}$ onto $\Ker T^*$ and $\Ker T$, respectively, it suffices to express bases for $\Ker T^*$ and $\Ker T$ in terms of the coefficient vector $\p=(p_0,p_1,\cdots,p_{n_1})^{\T}$ of the first characteristic polynomial.

\def\m{\mathbf{m}}
\def\n{\mathbf{n}}
Let
\begin{align*}
    \m=&(\overline{p_{n_1}},\overline{p_{n_1-1}},\cdots,\overline{p_0},\underbrace{0,\cdots,0}_{m-n_1-1})^{\T}\in\C^m, \\
    \n=&(p_0,p_1,\cdots,p_{n_1},\underbrace{0,\cdots,0}_{n-n_1-1})^{\T}\in\C^n.
\end{align*}
Let $J_m(0)$ and $J_n(0)$ be the nilpotent Jordan blocks of orders $m$ and $n$, respectively, that is,
\begin{align*}
    J_m(0)=\begin{pmatrix}
        0\\
        1&0\\
         &\ddots&\ddots\\
         &&1&0
    \end{pmatrix}\in\C^{m\times m},  J_n(0)=\begin{pmatrix}
        0\\
        1&0\\
         &\ddots&\ddots\\
         &&1&0
    \end{pmatrix}\in\C^{n\times n}. 
\end{align*}

Since $\dim \Ker T=n-\rank{T}=n-n_1$ and $\p\in\Ker T_{n_2-1,n_1+1}$, we have
\begin{align*}
    T\n=T(J_n(0)\n)=T(J_n(0)^2\n)=\cdots=T(J_n(0)^{n-n_1-1}\n)=0.
\end{align*}
It follows that the vectors $\n$, $ J_n(0)\n$, $ J_n(0)^2\n$, $\cdots$, $ J_n(0)^{n-n_1-1}\n$ form a basis for $\Ker T$.
Similarly, the vectors $\m$, $ J_m(0)\m$, $ J_m(0)^2\m$, $\cdots$, $ J_m(0)^{m-n_1-1}\m$ form a basis for $\Ker T^*$.

\section{Numerical examples}\label{Section4}

We first present an algorithm based on Theorem~\ref{thm.MPinv} for computing the Moore-Penrose inverse of an arbitrary $m\times n$ rank-deficient Toeplitz matrix; see Algorithm~\ref{alg:toeplitz-mpi}. Following Theorem~\ref{thm.MPinv}, we use the coefficient vector of the first characteristic polynomial to construct a quasi-direct decomposition of $T$ and the orthogonal projection matrices $\widetilde{M}$ and $\widetilde{N}$ onto $\Ker T^*$ and $\Ker T$, respectively, and thereby obtain the Moore-Penrose inverse of the Toeplitz matrix.

\begin{algorithm}[htbp]
\caption{Computing the Moore-Penrose inverse of a Toeplitz matrix via a quasi-direct decomposition}
\label{alg:toeplitz-mpi}
\KwIn{A nonzero rank-deficient Toeplitz matrix $T\in\mathbb C^{m\times n}$}
\KwOut{A numerical approximation $X\approx T^\dagger$ to the Moore-Penrose inverse}


$n_1=\rank T$, $\qquad n_2=m+n-n_1$


Choose $\p \in\Ker T_{n_2-1,n_1+1}\setminus\{0\}$

   $t=\min\{j:p(j)\ne0\}-1$, $\qquad s=\max\{j:p(j)\ne0\}-1$
   
   $\p_0= \p(t+1:s+1)$, $\qquad \p_0^f=flip(\p_0)$
   
   
 
   $P\gets I_n$, $\qquad Q\gets I_m$


   $P(i:i+s-t,i)\gets\p_0$, $\ i=1,\ldots,t$,   
   $\qquad P(i-(s-t):i,i)\gets \p_0$, $\ i=s+1,\ldots,n$
   
   $Q(i:i+s-t,i)\gets\p_0^f$, $\ i=1,\ldots,n_1-s$   
   $\qquad Q(i-(s-t):i,i)\gets\p_0^f$, $\ i=n_1-t+1,\ldots,m$

\BlankLine 
\tcp{Compute $A^{-1}$, $B^{-1}$, $C^{-1}$, and a $\{1,2\}$-inverse $G$ of $T$}  

$G\gets O_{n\times m}$

 \If{$s>t$}{$A= T(n_1-s+1:n_1-t,\ t+1:s)$
 
          $G(t+1:s,n_1-s+1:n_1-t) = A^{-1}$}
  
 \If{$n_1>s$}
       { $\beta=\mathcal L_{n_1-s}(\p_0^f)J_{n_1-s}T(1:n_1-s,\ n+t-s:n)\p_0$
       
       $B = \mathcal L_{n_1-s}(\beta)^{\T}$
       
       $G(n-n_1+s+1:n,1:n_1-s) = B^{-1}$ }

 \If{$t>0$}{ $\gamma=\mathcal L_{t}(\p_0)T(m-t+1:m,\ 1:s-t+1)\p_0$
 
 $C = \mathcal L_{t}(\gamma)$
 
         $G(1:t,m-t+1:m) = C^{-1}$}

   $G= PGQ^{\T}$

 \BlankLine
\tcp{Obtain the Moore-Penrose inverse }
  
  $MN \gets O_{(\max\{m,n\}-n_1+s-t)\times(\max\{m,n\}-n_1)}$
  
\For{$i = 1:\max\{m,n\}-n_1$}{$MN(i:s-t+i,i) = \p_0$}

$[U,~] = qr(MN,0)$

$N = [O_{t\times (n-n_1)};U(1:s-t+n-n_1,1:n-n_1);O_{(n_1-s)\times(n-n_1)}]$

$M = [O_{(n_1-s)\times(m-n_1)};flip(conj(U(1:s-t+m-n_1,1:m-n_1)));O_{t\times(m-n_1)}]$

$X \gets (I_n- NN^* )G(I_m - MM^*)$
   
\Return $X$

\begin{minipage}{\linewidth}
\footnotesize
\textbf{Notation: }$J_k$ is the flip matrix of order $k$. The notation $\mathcal L_k(v)$ denotes the lower triangular Toeplitz matrix of order $k$ whose first column is formed from the $l$-dimensional vector $v$ as follows: if $l\ge k$, the first column of $\mathcal L_k(v)$ consists of the first $k$ entries of $v$; if $l<k$, append $k-l$ entries equal to $0$ to $v$ to obtain a vector of length $k$, and use this vector as the first column.

\end{minipage}
\end{algorithm}

\subsection{Test matrices and experimental setup}

We use the rank-deficient Toeplitz matrices from Section~5 of~\cite{XW2016} to compare the numerical performance of the formula proposed here with that of the Gohberg-Semencul-Trench-type formula in that reference. The experiment examines computational error and runtime as the matrix order increases. We denote the result computed using Theorem~\ref{thm.MPinv} by $X_{\mathrm{QD}}$ and the result obtained using the comparison method by $X_{\mathrm{F}{\scriptscriptstyle\mathrm{M}}\mathrm{MP}}$. Both methods are rerun under the same hardware and software conditions.

Consider the Toeplitz matrix $T$ of order $n$ with first column and first row given, respectively, by
\begin{align*}
    \mathbf c_n&=\left(1,\frac12,\ldots,\frac1{n-1},1\right)^{\T},\\ 
    \mathbf r_n&=\left(1,\frac1{n-1},\frac1{n-2},\ldots,\frac12,1\right)
\end{align*}
The last row of this matrix is identical to its first row, and its last column is identical to its first column. 

Let $C$ be the leading principal submatrix of $T$ of order $n-1$. Then $C$ is an invertible circulant matrix with first column $(1,1/2,\ldots,1/(n-1))^{\T}$. Set
\[
    E_n=\begin{pmatrix}I_{n-1}\\e_1^{\T}\end{pmatrix},
    \qquad D_{n-1}=E_n^{\T}E_n=I_{n-1}+e_1e_1^{\T},
\]
where $e_1$ is the first standard basis vector in $\C^{n-1}$. Then
\begin{equation*}
    T=E_nCE_n^{\T},
\end{equation*}
and its Moore-Penrose inverse admits the following formula~\cite{XW2016}:
\begin{equation}\label{eq:num-reference}
    T^{\dag}=E_nD_{n-1}^{-1}C^{-1}D_{n-1}^{-1}E_n^{\T}.
\end{equation}
Formula~\eqref{eq:num-reference} holds in exact arithmetic. In the experiment, $C^{-1}$ is computed using the Fourier diagonalization of the circulant matrix to obtain a double-precision reference solution $X_{\mathrm{ref}}$ for $T^{\dag}$, and all Fourier eigenvalues are checked to be nonzero. This reference solution is used only for error assessment and is not used in the computations performed by either method.

The proposed method is based on Theorem~\ref{thm.MPinv}. We use the coefficient vector of the first characteristic polynomial to construct a $\{1,2\}$-inverse of $T$ and the orthogonal projection matrices $\widetilde{M}$ and $\widetilde{N}$ onto the null spaces $\Ker T^*$ and $\Ker T$, respectively.
For the orthogonal projections, we first compute orthonormal bases for $\Ker T^*$ and $\Ker T$, and form the matrices $M$ and $N$ with these basis vectors as columns, respectively. The orthogonal projection matrices are then given by $\widetilde{M}=MM^*$ and $\widetilde{N}=NN^*$.
A $\{1,2\}$-inverse $G$ of $T$ is computed from the coefficient vector of the first characteristic polynomial of $T$ using Theorem~\ref{thm.T12}.
Finally, the Moore-Penrose inverse of the Toeplitz matrix is computed using Theorem~\ref{thm.MPinv}.

The comparison method is the algorithm in~\cite{XW2016} for computing the Moore-Penrose inverse of a Toeplitz matrix using a Gohberg-Semencul-Trench-type formula, referred to as the F{\footnotesize M}MP algorithm.
Since the leading principal submatrix $C$ of $T$ of order $n-1$ is an invertible circulant matrix, the parameter required by the algorithm is $\alpha=1$ and need not be determined using LSQR. The implementation therefore requires four LSQR solves. The LSQR tolerance is set to $1\times 10^{-12}$ throughout, and the maximum number of iterations is set to the number of columns of the corresponding coefficient matrix.
All Toeplitz matrix-vector products in both methods are implemented using the fast Fourier transform.

All computations are performed in double-precision floating-point arithmetic using MATLAB R2024a on a 64-bit Windows system equipped with an Intel Core i5-10210U processor. All solves pass the convergence checks. Both implementations are warmed up before timing. For each matrix size, each method is run five times, and the median runtime is reported. The timings cover only the computation of the Moore-Penrose inverse and exclude the construction of the test matrix and the reference solution, as well as the error evaluation. Both methods explicitly output the entire $n\times n$ inverse matrix, so both require $\mathcal{O}(n^2)$ storage for the output.

Computational accuracy is measured by the absolute error in the matrix 2-norm,
\begin{equation*}
    E_2(X)=\|X-X_{\mathrm{ref}}\|_2 
\end{equation*}
where $X$ is either $X_{\mathrm{QD}}$ or $X_{\mathrm{F}{\scriptscriptstyle\mathrm{M}}\mathrm{MP}}$. Here, the matrix 2-norm is computed using MATLAB's \texttt{norm(X-X\textsubscript{ref},2)}, whereas~\cite{XW2016} estimates the corresponding error using \texttt{normest}; the resulting numerical values are therefore not expected to agree exactly.

\subsection{Runtime and error comparison for different matrix orders}

For $n=512$, $1024$, $2048$, $4096$, and $8192$, the absolute errors in the matrix 2-norm and the runtimes are reported in Table~\ref{tab:num-size}.
\begin{table}[htbp]
    \centering
    \caption{Errors and runtimes of the two methods. Times are in seconds and are the medians of five runs.}
    \label{tab:num-size}
    \begin{tabular}{cc cc c}
        \hline
        & \multicolumn{2}{c}{Proposed method} & \multicolumn{2}{c}{F{\footnotesize M}MP method}\\
        $n$ & Time & $E_2(X_{\mathrm{QD}})$ & Time & $E_2(X_{\mathrm{F}{\scriptscriptstyle\mathrm{M}}\mathrm{MP}})$\\
        \hline
        512  & 0.0306 & 3.0580E-15 & 0.0313   &3.2563E-11\\
        1024 & 0.0769 &2.0614E-15 &0.0562   &3.5965E-11 \\
        2048 & 0.2494 &6.4051E-15 & 0.1931 & 4.6960E-11\\
        4096 & 0.8924 &7.9610E-15 &0.4911 & 4.2858E-11\\
        8192 & 3.5909 & 6.6018E-15  &2.4575 & 8.5301E-11\\
        \hline
    \end{tabular}
\end{table}

For these matrix sizes, the absolute errors of the proposed method range from $2.06\times10^{-15}$ to $7.97\times10^{-15}$, whereas those of the F{\footnotesize M}MP method range from $3.25\times10^{-11}$ to $8.54\times10^{-11}$. Comparing the errors at each size shows that the errors of the proposed method are approximately $5.73\times10^{-5}$ to $1.86\times10^{-4}$ times those of the F{\footnotesize M}MP method, corresponding to a reduction of about four orders of magnitude. These results indicate that, for the present matrix family, the implementation of the proposed formula achieves smaller forward errors. Across the tested range of matrix orders, the errors remain on the order of $10^{-15}$ and show no sustained increase as the matrix order grows. 

In terms of runtime, the proposed method is slightly faster than the F{\footnotesize M}MP method when $n=512$, taking $0.0306$ seconds and $0.0313$ seconds, respectively. For the other four matrix sizes, the proposed method takes approximately $1.29$ to $1.82$ times as long as the F{\footnotesize M}MP method. In particular, the runtimes are $0.8924$ seconds and $0.4911$ seconds for $n=4096$, and $3.5909$ seconds and $2.4575$ seconds for $n=8192$, respectively. The current implementation of the proposed algorithm still requires Toeplitz block inversion, matrix products, and orthogonal projections, so a shorter overall runtime is not guaranteed directly by the formula. All reported runtimes are medians of five runs and reflect the performance of the current implementations when producing the entire inverse matrix.

\section{Conclusions}\label{Section5}
This work focuses on constructing a new formula for the Moore-Penrose inverse of a rank-deficient Toeplitz matrix, based on the theory of rational generating functions for Toeplitz matrices. 
We first investigate the quasi-direct decomposition of rank-deficient Toeplitz matrices. Using this decomposition, we then derive a new formula for the Moore-Penrose inverses of rank-deficient Toeplitz matrices. 
This formula applies to all rank-deficient Toeplitz matrices over the complex field, and requires only the coefficient vector of the first characteristic polynomial and some entries of the Toeplitz matrix itself. 
For the rank-deficient Toeplitz matrices from Section 5 in \cite{XW2016}, our algorithm based on Theorem\ref{thm.MPinv} requires slightly more computation time for large-scale matrices, but can reduce the numerical error by approximately four orders of magnitude, compared with the FMMP algorithm presented in the reference \cite{XW2016}.

\section*{Acknowledgments}
\addcontentsline{toc}{section}{Acknowledgments}

\biboptions{sort&compress}


\begin{thebibliography}{00}
\addcontentsline{toc}{section}{Reference}


\bibitem{BZH2026}B.H. Ban, X.Z. Zhan, Y.J. Hu, On the rational generating functions of Toeplitz matrices, Linear Algebra Appl., 733 (2026) 155-170. 


\bibitem{BCB2003} D.A. Bini, G. Codevico, M. Van Barel, Solving Toeplitz least squares problems by means of Newton’s iteration, Numerical Algorithms, 33 (2003) 93–103.



\bibitem{FB1999} J.A. Fessler, S.D. Booth, Conjugate-gradient preconditioning methods for shift-variant PET image reconstruction, IEEE Transactions on Image Processing, 8 (5) (1999) 688-699. 

\bibitem{F1984} M. Fiedler, Quasidirect decompositions of Hankel and Toeplitz matrices, Linear Algebra Appl. 61 (1984) 155-174.



\bibitem{GS1972} I. Gohberg, A. Semencul, On the inversion of finite Toeplitz matrices and their continuous analogs, Math. Issled. 7 (2) (1972) 201–223 (in Russian). 




\bibitem{H2004} G. Heinig, Fast Algorithms for Toeplitz Least Squares Problems, in: Current Trends in Operator Theory and its Applications, volume 149 of Oper. Theory Adv. Appl.,  Birkh\"{a}user, Basel, 2004, pp.167-197.

\bibitem{HH1994} G. Heinig, F. Hellinger, Moore–Penrose inversion of square Toeplitz matrices, SIAM J. Matrix Anal. Appl., 15(2) (1994) 418–450.



\bibitem{HJ1986} G. Heinig, U. Jungnickel, Hankel Matrices Generated by Markov Parameters, Hankel Matrix Extension, Partial Realization, and Pad\'{e}-Approximation, in: Operator Theory: Advances and Applications, vol. 19, Birkh\"{a}user, Basel, 1986, pp. 231-253.



\bibitem{HJ1986_1} G. Heinig, U. Jungnickel, Hankel matrices generated by the Markov parameters of rational functions, Linear Algebra Appl. 76 (1986) 121-135.

\bibitem{HR1984}  G. Heinig, K. Rost, Algebraic methods for Toeplitz-like matrices and operators, volume 13 of Oper. Theory Adv. Appl., Birkh\"{a}user Verlag, Basel, 1984. 

\bibitem{HR2010}  G. Heinig, K. Rost, Introduction to Bezoutians, in: D.A. Bini, V. Mehrmann, V. Olshevsky, E.E. Tyrtyshnikov, M. van Barel(Eds.), Numerical Methods for Structured Matrices and Applications, volume 199 of Oper. Theory Adv. Appl., Birkhäuser, Basel, 2010, pp. 25-118.



\bibitem{NP2014} M.K. Ng, J.Y. Pan, Weighted Toeplitz regularized least squares computation for image restoration, SIAM J. Sci. Comput. 36 (1) (2014) B94-B121. 



\bibitem{WCN2004} Y.M. Wei, J.F. Cai, M. K. Ng, Computing Moore-Penrose inverses of Toeplitz matrices by Newton’s iteration, Mathematical and Computer Modelling, 40(1-2) (2004) 181–191.

\bibitem{XW2016} P.P. Xie, Y.M. Wei, The stability of formulae of the Gohberg–Semencul–Trench type for Moore–Penrose and group inverses of Toeplitz matrices, Linear Algebra Appl., 498 (2016) 117–135.

\bibitem{X1992} Z. Xu, On Moore-Penrose inverses of Toeplitz matrices, Linear Algebra Appl., 169 (1992) 9–15.

\end{thebibliography}
\end{document}